\documentclass{article}

\def\ra{\rightarrow}

\def\ss{\subseteq}

\def\e{\epsilon}

 \def\HollowBox #1#2{{\dimen0=#1 \advance\dimen0 by -#2       
       \dimen1=#1 \advance\dimen1 by #2                       
        \vrule height #1 depth #2 width #2                    
        \vrule height 0pt depth #2 width #1                   
        \llap{\vrule height #1 depth -\dimen0 width \dimen1}%
       \hskip -#2                                             
       \vrule height #1 depth #2 width #2}}                   
 \def\BoxOpTwo{\mathord{\HollowBox{6pt}{.4pt}}\;}             

\def\endpf{\hfill $\BoxOpTwo$}

\def\bomega{\partial \Omega}
\def\sss{\subset \, \, \subset}

\font\teneufm=eufm10
\font\seveneufm=eufm7
\font\fiveeufm=eufm5
\newfam\eufmfam
\textfont\eufmfam=\teneufm
\scriptfont\eufmfam=\seveneufm
\scriptscriptfont\eufmfam=\fiveeufm

\newfam\msbfam
\font\tenmsb=msbm10  \textfont\msbfam=\tenmsb
\font\sevenmsb=msbm7  \scriptfont\msbfam=\sevenmsb
\font\fivemsb=msbm5    \scriptscriptfont\msbfam=\fivemsb
\def\Bbb{\fam\msbfam \tenmsb}

\def\O{\Omega}

\def\RR{{\Bbb R}}
\def\CC{{\Bbb C}}

\def\NN{{\Bbb N}}

\newtheorem{theorem}{Theorem}
\newtheorem{corollary}[theorem]{Corollary}
\newtheorem{proposition}[theorem]{Proposition}
\newtheorem{lemma}[theorem]{Lemma}

\newtheorem{definition}{Definition}
\newtheorem{example}[definition]{EXAMPLE}

\makeindex

\usepackage{graphicx}

\usepackage{amsmath}

\begin{document}

\begin{center}
\huge \bf
Convexity in Real Analysis\footnote{{\bf MR Classfication Numbers:}
26B25, 52A05, 26B10, 26B35.} \end{center}
\vspace*{.12in}

\begin{center}
\large Steven G. Krantz\footnote{Author supported in part
by the National Science Foundation and by the Dean of the Graduate
School at Washington University.}
\end{center}
\vspace*{.15in}

\begin{center}
\today
\end{center}
\vspace*{.2in}

\begin{quotation}
{\bf Abstract:} \sl
We treat the classical notion of convexity in the context of
hard real analysis.  Definitions of the concept are given
in terms of defining functions and quadratic forms, and characterizations are provided
of different concrete notions of convexity.   This analytic
notion of convexity is related to more classical geometric ideas.
Applications are given both to analysis and geometry.
\end{quotation}
\vspace*{.25in}

\setcounter{section}{-1}

\section{Introduction}

Convexity is an old subject in mathematics.  Archimedes used convexity in
his studies of area and arc length.  The concept appeared intermittently
in the work of Fermat, Cauchy, Minkowski, and others.  Even Johannes Kepler
treated convexity.   But it can be said that the subject was not
really formalized until the seminal tract of Bonneson and Fenchel [BOF].
See also [FEN] for the history.	 Modern treatments of convexity
may be found in [LAY] and [VAL].

In what follows, we let the term ``domain'' denote a connected, open set.
We usually denote a domain by $\Omega$.
If $\Omega$ is a domain and $P, Q \in \Omega$ then the {\it closed segment}
determined by $P$ and $Q$ is the set
$$
\overline{PQ} \equiv \{(1 - t)P + tQ: 0 \leq t \leq 1\} \, .
$$

Most of the classical treatments of convexity rely on the following synthetic
definition of the concept:

\begin{definition} \rm
Let $\Omega \ss \RR^N$ be a domain.  We say that $\Omega$ is {\it convex}
if, whenever $P, Q \in \Omega$, then the closed segment $\overline{PQ}$ from $P$ to $Q$
lies in $\Omega$.
\end{definition}

Works such as [LAY] and [VAL] treat theorems of Helly and
Kirchberger---about configurations of convex sets in the plane, and points
in those convex sets. However, studies in analysis and differential
geometry (as opposed to synthetic geometry) require results---and
definitions---of a different type. We need hard analytic facts about the
shape of the boundary---formulated in differential-geometric language. We
need invariants that we can calculate and estimate. That is the point of
view that we wish to explore in the present paper.

\section{The Concept of Defining Function}

Let $\Omega \ss \RR^N$ be a domain with $C^1$ boundary.  A $C^1$ function
$\rho: \RR^N \ra \RR$ is called a {\it defining function}
for $\Omega$ if

\begin{enumerate}
\item $\Omega = \{ x \in \RR^N: \rho(x) < 0\}; $
\item $\mbox{}^c\overline{\Omega} = \{x \in \RR^N: \rho(x) > 0\}; $
\item $\nabla \rho(x) \not = 0\ \ \forall x \in \bomega$.
\end{enumerate}
In case $k \geq 2$ and $\rho$ is $C^k$ then we say
that the domain $\Omega$ has $C^k$ boundary.

This last point merits some discussion. For the notion of a domain having
$C^k$ boundary has many different formulations. One may say that $\Omega$
has $C^k$ boundary if $\partial \Omega$ is a regularly imbedded $C^k$
manifold in $\RR^N$. Or if $\partial \Omega$ is locally the graph of a
$C^k$ function. In the very classical setting of $\RR^2$, it is common to
say that the boundary of a domain or region (which of course is simply a
{\it curve} $\gamma: S^1 \ra \RR^2$) is $C^k$ if {\bf (a)} $\gamma$ is a
$C^k$ function and {\bf (b)} $\gamma' \ne
0$.

We shall not take the time here to prove the equivalence of all the
different formulations of $C^k$ boundary for a domain (but
see the rather thorough discussion in Appendix I of [KRA1]).  But
we do discuss the equivalence of the ``local graph'' definition
with the defining function definition.

First suppose that $\Omega$ is a domain with $C^k$ defining function $\rho$
as specified above, and let $P \in \partial \Omega$.  Since $\nabla \rho(P) \ne 0$,
the implicit function theorem (see [KRP2]) guarantees that there is a 
a neighborhood $V_P$ of $P$, a variable
(which we may take to be $x_N$) and a $C^k$ function $\varphi_P$ defined
on a small open set $U_P \ss \RR^{N-1}$ so that
$$
\partial \Omega \cap V_P = \{(x_1, x_2, \dots, x_N): x_N = 
\varphi_P(x_1, \dots, x_{N-1}) \ , \quad (x_1, \dots, x_{N-1}) \in U_P\} \, .
$$ 
Thus $\partial \Omega$ is locally the graph of the function $\varphi_P$ near $P$.

Conversely, assume that each point $P \in \partial \Omega$ has a neighborhood $V_P$
and an associated $U_P \ss \RR^{N-1}$ and function $\varphi_P$ such that
$$
\partial \Omega \cap V_P = \{(x_1, x_2, \dots, x_N): x_N = 
\varphi_P(x_1, \dots, x_{N-1}) \ , \quad (x_1, \dots, x_{N-1}) \in U_P\} \, .
$$	   
We may suppose that the positive $x_N$-axis points {\it out} of the domain,
and set $\rho_P(x) = x_N - \varphi_P(x_1,\dots, x_{N-1})$.  Thus, on
a small neighborhood of $P$, $\rho_P$ behaves like a defining function.
It is equal to 0 on the boundary, certainly has non-vanishing gradient,
and is $C^k$.

Now $\partial \Omega$ is compact, so we may cover $\partial \Omega$ with
finitely many $V_{P_1}, \dots, V_{P_k}$.  Let $\{\psi_j\}$ be a partition
of unity subordinate to this finite cover, and set
$$
\widetilde{\rho}(x) = \sum_{j=1}^k \psi_j(x) \cdot \rho_{P_j}(x) \, .
$$
Then, in a neighborhood of $\partial \Omega$, $\widetilde{\rho}$ is a defining function.
We may extend $\widetilde{\rho}$ to all of space as follows.  Let $V$ be a neighborhood
of $\partial \Omega$ on which $\widetilde{\rho}$ is defined.  Let $V'$ be an open, relatively
compact subset of $\Omega$ and $V''$ an open subset of ${}^c \overline{\Omega}$ 
so that $V, V', V''$ cover $\CC^n$.  Let $\eta, \eta', \eta''$ be a partition
of unity subordinate to the cover $V, V', V''$. Now set
$$
\rho(x) = \eta'(x) \cdot [-(C + 10)^2] + \eta(x) \cdot \widetilde{\rho}(x) + \eta''(x) \cdot  (C + 10)^2 \, .
$$
Here $C$ is a large positive constant that exceeds the diameter of $\Omega$.  Then $\rho$ is a globally
defined, $C^k$ function that is a defining function for $\Omega$.

\begin{definition}	\rm
Let $\Omega \ss \RR^N$ have $C^1$ boundary and let $\rho$
be a $C^1$ defining function.  Let $P \in \partial \Omega.$
An $N-$tuple $w = (w_1,\dots,w_N)$ of real numbers 
is called a {\it tangent vector} 
to $ \partial \Omega$ at $P$ if 
$$
 \sum_{j=1}^N (\partial \rho/\partial x_j)(P) \cdot w_j = 0. 
$$
We write $w \in T_P(\partial \Omega).$
\end{definition}

For $\Omega$ with $C^1$ boundary, we think of 
$\nu_P = \nu = \langle \partial \rho/\partial x_1(P), \dots, \partial \rho/\partial x_N(P) \rangle$ as
the outward-pointing normal vector to $\partial \Omega$ at $P$.
Of course the union of all the tangent vectors to $\partial\Omega$ at a point $P \in \partial \Omega$ 
is the {\it tangent plane} or {\it tangent hyperplane}.  The tangent hyperplane is defined by the condition
$$
\nu_P \cdot w = 0 \, .
$$
This definition makes sense when $\nu_P$ is well defined, in particular when $\partial \Omega$ is $C^1$.  

If $\Omega$ is convex and $\partial \Omega$ is not smooth---say
that it is Lipschitz---then any point $P \in \partial \Omega$
will still have one (or many) hyperplanes ${\cal P}$ such that
${\cal P} \, \cap \, \overline{\Omega} = \{P\}$. We call such
a hyperplane a {\it support hyperplane} for $\partial \Omega$
at $P$. As noted, such a support hyperplane need not be
unique. For example, if $\Omega = \{(x_1, x_2): |x_1| < 1,
|x_2| < 1\}$ then the points of the form $(\pm 1, \pm 1)$ in
the boundary do {\it not} have well-defined tangent planes,
but they do have (uncountably) many support hyperplanes.

Of course the definition of the normal $\nu_P$ makes sense only if it is independent
of the choice of $\rho.$  We shall address that issue in a moment.
It should be observed that the condition defining tangent vectors simply
mandates that $w \perp \nu_P$ at $P.$   And, after all, we know from calculus
that $\nabla \rho$ is the normal $\nu_P$ and that the normal is uniquely
determined and independent of the choice of $\rho.$  In principle, this
settles the well-definedness issue.

However this point is so important, and the point of view that
we are considering so pervasive, that further discussion is warranted.
The issue is this:  if $\widehat{\rho}$ is another defining function for
$\Omega$ then it should give the same tangent vectors as $\rho$ at any
point $P \in \bomega.$  The key to seeing that this is so is to
write $\widehat{\rho}(x) = h(x) \cdot \rho(x),$ for $h$ a function
that is non-vanishing near $\bomega.$  Then, for $P \in \bomega,$
\begin{align}
  \sum_{j=1}^N (\partial \widehat{\rho}/\partial x_j)(P) \cdot w_j
     & =  h(P) \cdot \left ( \sum_{j=1}^N (\partial \rho/\partial x_j)(P) \cdot w_j \right ) \notag \\
     &    \mbox{} + \rho(P) \cdot \left ( \sum_{j=1}^N (\partial h/\partial x_j)(P) \cdot w_j \right ) \notag  \\
     & =  h(P) \cdot \left ( \sum_{j=1}^N (\partial \rho/\partial x_j)(P) \cdot w_j \right )  \notag \\
     &    \qquad   + 0 , \tag{1.1} \\ \notag
\end{align}
\vspace*{.12in}

\noindent because $\rho(P) = 0.$   Thus $w$ is a tangent vector at $P$ vis a vis $\rho$
if and only if $w$ is a tangent vector vis a vis $\widehat{\rho}.$  But
why does $h$ exist?

After a change of coordinates, it is enough to assume that we are
dealing with a piece of $\bomega$ that is a piece of flat, $(N-1)-$dimensional
real hypersurface (just use the implicit function theorem).  Thus we
may take $\rho(x) = x_N$ and $P = 0.$  Then any other defining function 
$\widehat{\rho}$ for $\bomega$ near $P$ must have the Taylor expansion
$$
  \widehat{\rho}(x) = c \cdot x_N + {\cal R}(x)  \eqno (1.2)
$$
about $0.$  Here ${\cal R}$ is a remainder term\footnote{We may think of (1.2)
as proved by integration by parts in the $x_N$ variable only, and that gives this favorable estimate
on the error terms ${\cal R}(x)$.} satisfying ${\cal R}(x) = o(|x_N|)$.
There is no loss of generality to take $c = 1,$ and we do so
in what follows.  Thus we wish to define 
$$
  h(x) = \frac{\widehat{\rho}(x)}{\rho(x)} = 1 + {\cal S}(x) . 
$$
Here ${\cal S}(x) \equiv {\cal R}(x)/x_N$ and ${\cal S}(x) = o(1)$
as $x_N \ra 0.$
Since this remainder term involves a derivative of $\widehat{\rho},$
it is plain that $h$ is not even differentiable.  (An explicit counterexample
is given by $\widehat{\rho}(x) = x_N \cdot(1 + |x_N|).$)  Thus the 
program that we attempted in equation $(1.1)$ above is apparently flawed.

\setcounter{theorem}{2}

However an inspection of the explicit form of the remainder term
${\cal R}$ reveals that, because $\widehat{\rho}$ is constant on $\bomega,$
$h$ as defined above {\it is} continuously 
differentiable {\it in tangential directions.}
That is, for tangent vectors $w$ (vectors that are orthogonal 
to $\nu_P$), the derivative
$$
 \sum_j \frac{\partial h}{\partial x_j} (P) w_j 
$$
{\it is} defined.  Thus it does indeed turn out that our definition of 
tangent vector is well-posed when it is applied to vectors {\it that are
already known to be tangent vectors} by the geometric definition
$w \cdot \nu_P = 0.$  
For vectors that are {\it not}
geometric tangent vectors, an even simpler argument shows that
$$
 \sum_j \frac{\partial \widehat{\rho}}{\partial x_j} (P) w_j \not = 0
$$
if and only if
$$
 \sum_j \frac{\partial \rho}{\partial x_j} (P) w_j \not = 0 .
$$
Thus Definition 2 is well-posed.
Questions similar to the one just discussed will come up below when
we define convexity using $C^2$ defining functions.  They are resolved
in just the same way and we shall leave details to the reader.

The reader should check that the discussion above proves the following:
if $\rho, \widetilde{\rho}$ are $C^k$ defining functions for a domain
$\Omega$, with $k \geq 2$, then there is a $C^{k-1}$, nonvanishing 
function $h$ defined near $\partial \Omega$
such that $\rho = h \cdot \widetilde{\rho}.$
\medskip \\

\section{The Analytic Definition of Convexity}

\indent For convenience, we restrict attention for
this section to {\it bounded}
domains.  Many of our definitions would need to be modified, and
extra arguments given in proofs, were we to consider unbounded
domains as well.

\begin{definition}	\rm
Let $\Omega \sss \RR^N$ be a domain with $C^2$ boundary and $\rho$ a
defining function for $\Omega.$  Fix a point $P \in \partial \Omega.$
We say that $\partial \Omega$ is analytically (weakly) {\it convex} at $P$ if
$$
 \sum_{j,k=1}^N \frac{\partial^2\rho}{\partial x_j \partial x_k}(P) w_j w_k \geq 0,\ \ 
                     \forall w \in T_P(\partial \Omega) . 
$$
We say that $\partial \Omega$ is analytically {\it strongly (strictly) convex} at $P$
if the inequality is strict whenever $w \not = 0.$

If $\partial \Omega$ is convex (resp. strongly convex) at each boundary 
point then we say that
$\Omega$ is convex (resp. strongly convex).  
\end{definition}  

One interesting and useful feature of this new definition of convexity is
that it treats the concept point-by-point.  The classical, synthetic definition
specifies convexity for the whole domain at once.

It is natural to ask whether the new definition of convexity is independent of
the choice of defining function.  We have the following result:

\begin{proposition} \sl
Let $\Omega \ss \RR^N$ be a domain with $C^2$ boundary.  Let $\rho$ and $\rho'$ be 
$C^2$ defining functions for $\Omega$, and assume that, at points $x$ near $\partial \Omega$,
$$
\rho(x) = h(x) \cdot \rho'(x)
$$
for some non-vanishing, $C^2$ function $h$.  Let $P \in \partial \Omega$.  Then $\Omega$
is convex at $P$ when measured with the defining function $\rho$ if and only
if $\Omega$ is convex at $P$ when measured with the defining function $\rho'$.
\end{proposition}
{\bf Proof:}  We calculate that
\begin{eqnarray*}
\frac{\partial^2}{\partial x_j \partial x_k} \rho(P) & = &                                               
	   h(P) \cdot \frac{\partial^2 \rho'}{\partial x_j \partial x_k} (P) + 
	   \rho'(P) \cdot \frac{\partial^2 h}{\partial x_j \partial x_k} (P) 	\\
	  && \qquad  + \frac{\partial \rho'}{\partial x_j} (P) \frac{\partial h}{\partial x_k} (P) + 
	   \frac{\partial \rho'}{\partial x_k} (P) \frac{\partial h}{\partial x_j} (P)  \\
	 & = &  h(P) \cdot \frac{\partial^2 \rho'}{\partial x_j \partial x_k} (P) + 
   	   \frac{\partial \rho'}{\partial x_j} (P) \frac{\partial h}{\partial x_k} (P) + 
	   \frac{\partial \rho'}{\partial x_k} (P) \frac{\partial h}{\partial x_j} (P)   \\
\end{eqnarray*}
because $\rho'(P) = 0$.  But then, if $w$ is a tangent vector to $\partial \Omega$ at $P$, we see that
\begin{eqnarray*}
\sum_{j, k} \frac{\partial^2}{\partial x_j \partial x_k} \rho(P) w_j w_k & = & 
         h(P) \sum_{j, k} \frac{\partial^2 \rho'}{\partial x_j \partial x_k} (P) w_j w_k \\
        && \qquad + \left [ \sum_{j} \frac{\partial \rho'}{\partial x_j} (P) w_j \right ] \left [ \sum_{k} \frac{\partial h}{\partial x_k} (P) w_k \right ] \\
                && \qquad + \left [ \sum_{k} \frac{\partial \rho'}{\partial x_k} (P) w_k \right ] \left [ \sum_{j} \frac{\partial h}{\partial x_j} (P) w_j \right ]  \, .
\end{eqnarray*}
If we suppose that $P$ is a point of convexity relative to the defining function $\rho'$, then the first sum is nonnegative.
Of course $h$ is positive, so the first expression is then $\geq 0$.  Since $w$ is a tangent vector, the sum in $j$ in the second expression
vanishes.  Likewise the sum in $k$ in the third expression vanishes.

In the end, we see that the Hessian of $\rho$ is positive semi-definite on the tangent space if the Hession of $\rho'$ is.
The reasoning also works if the roles of $\rho$ and $\rho'$ are reversed.  The result is thus proved.
\endpf 
\smallskip \\

The quadratic form
$$
 \left ( \frac{\partial^2\rho}{\partial x_j \partial x_k}(P) \right )_{j,k=1}^N 
$$
is frequently called the ``real Hessian'' of the function $\rho.$  This form
carries considerable geometric information about the boundary
of $\Omega.$  It is of course closely related to the second fundamental
form of Riemannian geometry (see B. O'Neill [ONE]). 

There is a technical difference between ``strong'' and ``strict'' convexity
that we shall not discuss here (see L. Lempert [LEM] for details).  It is common
to use either of the words ``strong'' or ``strict'' to mean that the 
inequality in the last definition is strict when $w \not = 0.$
The reader may wish to verify for himself that, at a strongly
convex boundary point, all curvatures are positive (in fact one
may, by the positive definiteness of the matrix 
$\left ( \partial^2 \rho/\partial x_j \partial x_k \right ),$
impose a change of coordinates at $P$ so that the boundary
of $\Omega$ agrees with a ball up to second order at $P$).  

Now we explore our analytic notions of convexity.  The first lemma is
a technical one:

\begin{lemma} \sl    
Let $\Omega \ss\RR^N$ be strongly convex.  Then there is a constant
$C > 0$ and a defining function
$\widetilde{\rho}$ for $\Omega$ such that
$$
  \sum_{j,k=1}^N  \frac{\partial^2\widetilde{\rho}}{\partial x_j \partial x_k} (P) w_j w_k \geq C |w|^2 , \ \ 
                                               \forall P \in \partial \Omega, w \in \RR^N . \eqno (2.1) 
$$
\end{lemma}
{\bf Proof:}  Let $\rho$ be some fixed ${C}^2$ defining function for $\Omega.$
For $\lambda > 0$ define
$$
  \rho_\lambda(x) = \frac{\mbox{\rm exp}(\lambda \rho(x)) - 1}{\lambda} . 
$$
We shall select $\lambda$ large in a moment.  Let $P \in \bomega$ and set
$$
 X = X_P = \left \{ w \in \RR^N: |w| = 1\ \mbox{\rm and}\ 
 \sum_{j,k} \frac{\partial^2 \rho}{\partial x_j \partial x_k}(P) w_j w_k \leq 0 \right \} . 
$$
Then no element of $X$ could be a tangent vector at $P,$ 
hence $X \ss \{w: |w| = 1\ \mbox{\rm and}\ \sum_j \partial \rho/\partial x_j(P) w_j \not = 0\}.$
Since $X$ is defined by a non-strict inequality, it is closed; it is
of course also bounded.  Hence $X$ is compact and
$$
 \mu \equiv \min \left \{\left | \sum_j \partial \rho/\partial x_j(P) w_j\right |: w \in X \right \} 
$$
is attained and is non-zero.  Define
$$
  \lambda = \frac{-\min_{w \in X}\sum_{j,k} \frac{\partial^2 \rho}{\partial x_j \partial x_k}(P)w_j w_k}{\mu^2}  + 1 . 
$$
Set $\widetilde{\rho} = \rho_\lambda.$  Then for any $w \in \RR^N$ with 
$|w| = 1$ we have (since $\mbox{\rm exp}(\rho_(P)) = 1$)
that
\begin{eqnarray*}
   \sum_{j,k} \frac{\partial^2 \widetilde{\rho}}{\partial x_j \partial x_k}(P)w_j w_k 
            & = & \sum_{j,k} \left \{ \frac{\partial^2 \rho}{\partial x_j \partial x_k}(P) + 
                \lambda \frac{\partial \rho}{\partial x_j}(P) \frac{\partial \rho}{\partial x_k}(P) \right \} w_j w_k \\
            & = & \sum_{j,k} \left \{ \frac{\partial^2 \rho}{\partial x_j \partial x_k} \right \} (P) w_j w_k + 
            \lambda \left | \sum_j \frac{\partial \rho}{\partial x_j} (P) w_j \right |^2   
\end{eqnarray*}
If $w \not \in X$ then this expression is positive by definition.
If $w \in X$ then the expression is positive by the choice of $\lambda.$
Since $\{w \in \RR^N: |w| = 1\}$ is compact, there is thus a $C > 0$
such that
$$
 \sum_{j,k} \left \{ \frac{\partial^2 \widetilde{\rho}}{\partial x_j \partial x_k} \right \} (P) w_j w_k \geq C, \ \ 
                                  \forall w \in \RR^N\ \mbox{\rm such that}\ |w| = 1 . 
$$

This establishes our inequality $(2.1)$ for $P \in \bomega$ fixed and $w$ in the
unit sphere of $\RR^N.$  For arbitrary $w,$ we set $w = |w| \widehat{w},$
with $\widehat{w}$ in the unit sphere.   Then $(2.1)$ holds for $\widehat{w}.$  Multiplying
both sides of the inequality for $\widehat{w}$ by $|w|^2$ and performing
some algebraic manipulations gives the result for fixed $P$ and
all $w \in \RR^N.$

Finally, notice that our estimates---in particular the existence of 
$C,$ hold uniformly over points in $\bomega$ near $P.$  Since $\bomega$
is compact, we see that the constant $C$ may be chosen uniformly over
all boundary points of $\Omega.$ \endpf 
\smallskip \\

Notice that the statement of the lemma has two important features:
{\bf (i)} that the constant $C$ may be selected uniformly over the boundary
and {\bf (ii)} that the inequality $(2.1)$ holds for all $w \in \RR^N$ (not
just tangent vectors).  In fact it is impossible to arrange 
for anything like $(2.1)$ to be true at a weakly convex point.  

Our proof shows in fact that $(2.1)$ is true not just for $P \in \bomega$
but for $P$ in a neighborhood of $\bomega.$  It is this sort of
stability of the notion of strong convexity that makes it a more
useful device than ordinary (weak) convexity.

\begin{proposition}    
If $\Omega$ is strongly convex then $\Omega$ is geometrically convex.
\end{proposition}
{\bf Proof:}  We use a connectedness argument.  

Clearly
$\Omega \times \Omega$ is connected.
Set $S = \{(P_1,P_2) \in \Omega \times \Omega: (1 - \lambda) P_1 + \lambda P_2 \in \Omega,
                      \ \mbox{\rm all}\ 0 < \lambda < 1\}.$  Then $S$ 
is plainly open and non-empty.  

To see that $S$ is closed, fix a defining function $\widetilde{\rho}$ for
$\Omega$ as in the Lemma.
If $S$ is not closed in $\Omega \times \Omega$ then there 
exist $P_1, P_2 \in \Omega$ such that the function
$$
 t \mapsto \widetilde{\rho}((1 - t) P_1 + t P_2) 
$$
assumes an interior maximum value of $0$ on $[0,1].$  But the positive definiteness
of the real Hessian of $\widetilde{\rho}$ contradicts that assertion.
The proof is complete. \endpf  \smallskip \\ 

We gave a special proof that strong convexity implies geometric
convexity simply to illustrate the utility of the strong
convexity concept. It is possible to prove that an arbitrary
(weakly) convex domain is geometrically convex by showing that
such a domain can be written as the increasing union of
strongly convex domains. However the proof is difficult and
technical. We thus give another proof of this fact:

\begin{proposition}   
If $\Omega$ is (weakly) convex then $\Omega$ is geometrically convex.
\end{proposition}
{\bf Proof:}  To simplify the proof we shall assume that $\Omega$
has at least ${C}^3$ boundary.

Assume without loss of generality that $N \geq 2$ and
$0 \in \O.$  For $\e > 0,$
let $\rho_\e(x) = \rho(x) + \e |x|^{2M}/M$ and $\O_\e = \{x: \rho_\e(x) < 0\}.$
Then $\O_\e \ss \O_{\epsilon'}$ if $\e' < \e$ and $\cup_{\e > 0} \O_\e = \O.$
If $M \in \NN$ is large and $\e$ is small, then $\O_\e$ is strongly
convex.  By Proposition 4, each ${\O}_\e$ is geometrically convex, so $\O$ is
convex. 
\endpf  \smallskip \\ 

\indent We mention in passing that a nice treatment of convexity, from roughly
the point of view presented here, appears in V. Vladimirov [VLA].

\begin{proposition}    
Let $\Omega \sss \RR^N$ have $C^2$ boundary and be geometrically convex.
Then $\Omega$ is (weakly) convex.
\end{proposition}
{\bf Proof:}  Seeking a contradiction, we suppose that for 
some $P \in \partial \Omega$ and some $w \in T_P(\partial \Omega)$
we have
$$
 \sum_{j,k} \frac{\partial^2 \rho}{\partial x_j \partial x_k} (P) w_j w_k = - 2K < 0 . \eqno (2.2) 
$$
Suppose without loss of generality that coordinates have been selected
in $\RR^N$ so that $P = 0$ and $(0,0, \dots,0,1)$ is the unit outward normal
vector to $\bomega$ at $P.$  We may further normalize the
defining function $\rho$ so that $\partial \rho/\partial x_N(0) = 1.$
Let $Q = Q^t = tw + \epsilon \cdot (0, 0,\dots, 0,1),$ where
$\epsilon > 0$ and $t \in \RR.$  Then, by Taylor's expansion,
\begin{eqnarray*}
   \rho(Q) & = & \rho(0) + \sum_{j=1}^N \frac{\partial \rho}{\partial x_j} (0) Q_j + 
                     \frac{1}{2} \sum_{j,k=1}^N  \frac{\partial^2 \rho}{\partial x_j \partial x_k} (0) Q_j Q_k 
                     + o(|Q|^2) \\
           & = & \epsilon \frac{\partial \rho}{\partial x_N}(0) +
                    \frac{t^2}{2} \sum_{j,k=1}^N  \frac{\partial^2 \rho}{\partial x_j \partial x_k} (0) w_j w_k 
                + {\cal O}(\epsilon^2) + o(t^2) \\
           & = & \epsilon - K t^2 + {\cal O}(\epsilon^2) + o(t^2) .
\end{eqnarray*}
Thus if $t = 0$ and $\epsilon > 0$ is small enough then $\rho(Q) > 0.$
However, for that same value of $\epsilon,$ if $|t| > \sqrt{2\epsilon/K}$
then $\rho(Q) < 0.$
This contradicts the
definition of geometric convexity.
\endpf  \smallskip \\

\paragraph{Remark:} The reader can already see in the proof of the proposition
how useful the 
quantitative version of convexity can be.  

The assumption that
$\partial \Omega$ be $C^2$ is not very restrictive, for convex
functions of one variable are twice differentiable almost everywhere
(see A. Zygmund [ZYG]).  On the other hand, $C^2$ smoothness of the
boundary is essential for our approach to the subject.
\endpf 
\medskip \\

\paragraph{Exercises for the Reader:}  

If $\Omega \ss \RR^N$ is a domain then
the {\it closed convex hull} of $\Omega$ is defined to be the closure
of the set $\bigl \{\sum_{j=1}^m \lambda_j s_j: s_j \in \Omega, m \in \NN, \lambda_j \geq 0, \sum \lambda_j = 1 \bigr \}.$
Equivalently, the closed convex hull of $\Omega$ is the intersection of all
closed, convex sets that contain $\Omega$.

Assume in the following problems 
that $\overline{\Omega} \ss \RR^N$ is closed, bounded, and convex.  Assume
that $\Omega$ has $C^2$ boundary.
\smallskip \\

\noindent ({\bf a}) \ \  We shall say more about extreme points in the penultimate section.
For now, a point $P \in \partial \Omega$ is extreme (for $\Omega$ convex) if,
whenever $P = (1 - \lambda)x + \lambda y$ and $0 \leq \lambda \leq 1$, $x, y \in \overline{\Omega}$,
then $x = y = P$.  Prove that $\overline{\Omega}$ is the closed convex hull
of its extreme points (this result is usually referred to as the
{\it Krein-Milman theorem} and is true in much greater generality).

\noindent ({\bf b}) \ \ Let $P \in \bomega$ be extreme.  
Let ${\bf p} = P + T_P(\bomega)$ be the geometric tangent affine
hyperplane to the boundary of $\Omega$ that passes through $P.$ 
Show by an example that it
is not necessarily the case that ${\bf  p} \cap \overline{\Omega} = \{P\}.$

\noindent ({\bf c}) \ \ Prove that if $\Omega_0$ is {\it any} bounded domain with
$C^2$ boundary then there is a relatively open subset $U$ of $\bomega_0$
such that $U$ is strongly convex.  (Hint:  Fix $x_0 \in \Omega_0$
and choose $P \in \bomega_0$ that is as far as possible from $x_0).$

\noindent ({\bf d})  \ \ If $\Omega$ is a convex domain then the Minkowski functional\footnote{A simple instance
of the Minkowski functional is this.  Let $K \ss \RR^N$ be convex.  For $x \in \RR^N$, define
$$
p(x) = \inf \{r > 0: x \in rK\} \, .
$$
Then $p$ is a Minkowski functional for $K$.}
(see [LAY]) less 1 gives a convex defining function for $\Omega.$
\medskip \\

\section{Convex Functions and Exhaustion Functions}

Let $F: \RR^N \ra \RR$ be a function.  We say that
$F$ is {\it convex} if, for any $P, Q \in \RR^N$ and
any $0 \leq t \leq 1$, it hold that
$$
F((1 - t) P + tQ) \leq (1 - t)f(P) + t f(Q) \, .
$$
In the case that $F$ is $C^2$, we may restrict $F$ to the 
line passing through $P$ and $Q$ and differentiate
the function 
$$
\varphi_{P,Q}: t \longmapsto F((1 - t)P + tQ)
$$
twice to see (from calculus---reference [BLK]) that 
$$
\frac{d^2}{dt^2} \varphi_{P,Q} \geq 0 \, .
$$
If we set $\alpha = Q - P = (\alpha_1, \alpha_2, \dots, \alpha_N)$, then
this last result may be written as
$$
\frac{\partial^2}{\partial \alpha^2} F \geq 0 \, .
$$
This in turn may be rewritten as
$$
\sum_{j,k} \frac{\partial^2}{\partial x_j \partial x_k} \alpha_j \alpha_k \geq 0 \, .
$$
In other words, the Hessian of $F$ is positive semi-definite.

In the case that a $C^2$ function $F$ has positive definite Hessian
at each point then we say that $F$ is {\it strictly convex} or {\it strongly convex}.

The reasoning in the penultimate paragraph can easily be reversed to see that the
following is true:

\begin{proposition} \sl
A $C^2$ function on $\RR$ is convex if and only if it has positive semi-definite
Hessian at each point of its domain.
\end{proposition}

Of course it is also useful to consider convex functions on a domain.  Certainly
we may say that $F: \Omega \ra \RR$ is convex (with $\Omega$ a convex domain) if 
$$
\sum_{j,k} \frac{\partial^2}{\partial x_j \partial x_k} (x) \alpha_j \alpha_k \geq 0
$$
for all $(\alpha_1, \dots, \alpha_N)$ and all points $x \in \Omega$.  Equivalently,
$F$ is convex on a convex domain $\Omega$ if, whenever $P, Q \in \Omega$ and $0 \leq \lambda \leq 1$
we have
$$
F((1 - t) P + tQ) \leq (1 - t)f(P) + t f(Q) \, .
$$

It is straightforward to prove that any convex function is continuous.  See [ZYG]
or [VLA, p.\ 85].  Other properties of convex functions are worth noting.  For example,
if $f: \RR^n \ra \RR$ is convex and $\varphi: \RR \ra \RR$ is convex and increasing then $\varphi \circ f$
is convex.  Certainly the sum of any two convex functions is convex.  If $\{f_\alpha\}_{\alpha \in A}$ is any
family of convex functions then
$$
f(x) \equiv \sup_{\alpha \in A} f_\alpha(x)
$$
is convex.  The proof of this latter assertion is straightforward:
If $P, Q$ lie in the common domain of the $f_\alpha$ and $0 \leq \lambda \leq 1$  and $\alpha \in A$ then
$$
f_\alpha((1- \lambda)P + \lambda Q) \leq (1 - \lambda)f_\alpha(P) + \lambda f_\alpha(Q) \, .
$$
Then certainly
$$
f_\alpha((1- \lambda)P + \lambda Q) \leq (1 - \lambda)f(P) + \lambda f(Q) \, .
$$
Now take the supremum over $\alpha$ on the lefthand side to obtain the result.

It is always useful to be able to characterize geometric properties of domains
in terms of functions.  For functions are more flexible objects than domains:  one
can do more with functions.  With this thought in mind we make the following definition:

\begin{definition} \rm
Let $\Omega \ss \RR^N$ be a bounded domain.  We call a function
$$
\lambda: \Omega \ra \RR
$$
an {\it exhaustion function} if, for each $c \in \RR$, the set
$$
\lambda^{-1}((-\infty, c]) = \{x \in \Omega: \lambda(x) \leq c\}
$$
is a compact subset of $\Omega$.
\end{definition}
The key idea here is that the function $\lambda$ is real-valued and blows up
at $\partial \Omega$.

\begin{theorem} \sl
A domain $\Omega \ss \RR^N$ is convex if and only if it
has a continuous, convex exhaustion function.
\end{theorem}
{\bf Proof:}  If $\Omega$ possesses such an exhaustion function $\lambda$, then
 the domains
$$
 \Omega_k \equiv \{x \in \Omega: \lambda < k\}
$$
are convex.  And $\Omega$ itself is the increasing union of the $\Omega_k$.  It
follows immediately, from the synthetic definition of convexity, that $\Omega$
is convex.

For the converse, observe that if $\Omega$ is convex and $P \in \partial \Omega$, then
the tangent hyperplane at $P$ has the form $a \cdot (x - P) = 0$.  Here $a$ is a Euclidean
unit vector.   It then follows that the quantity $a \cdot (x - P)$ is the distance from
$x \in \Omega$ to this hyperplane.  Now the function
$$
\mu_{a,P}(x) \equiv - \ln a \cdot (x - P)
$$
is convex since one may calculate the Hessian ${\cal H}$ directly.  Its value at a point $x$ equals
$$
{\cal H}(b,b) = \frac{(a \cdot b)^2}{[a \cdot (x - P)]^2} \geq 0 \, .
$$

If $\delta_\Omega(x)$ is the Euclidean 
distance of $x$ to $\partial \Omega$, then
$$
- \log \delta_\Omega (x) = \sup_{P \in \partial \Omega} \left [ - \log a \cdot (x - P) \right ] \, .
$$
Thus $- \log \delta_\Omega$ is a convex function that blows up at $\partial \Omega$.  Now set
$$
\lambda(x) = \max \biggl \{ \tau_\Omega(x), |x|^2 \biggr \} \, .
$$
This is a continuous, convex function that blows up at the boundary.
So it is the convex exhaustion function that we seek.
\endpf 
\smallskip \\

\begin{lemma} \sl
Let $F$ be a convex function on $\RR^N$.  Then there is a sequences
$f_1 \geq f_2 \geq \cdots$ of $C^\infty$, strongly convex functions
such that $f_j \ra F$ pointwise.
\end{lemma}
{\bf Proof:}  Let $\varphi$ be a $C_c^\infty$ function which
is nonnegative and has integral 1.  We may also take $\varphi$ to
be supported in the unit ball, and be radial.  For $\epsilon > 0$
we set
$$
\varphi_\epsilon (X) = \epsilon^{-N} \varphi(x/\epsilon) \, .
$$
We define
$$
F_\epsilon(x) = F * \varphi_\epsilon (x) = \int F(x - t) \varphi_\epsilon (t) \, dt \, .
$$

We assert that each $F_\epsilon$ is convex.  For let $P, Q \in \RR^N$
and $0 \leq \lambda \leq 1$.  Then
\begin{eqnarray*} 
  F_\epsilon((1 - \lambda)P + \lambda Q) & = & \int F((1 - \lambda)P + \lambda Q - t) \varphi_\epsilon (t) \, dt \\
					 & = & \int F((1 - \lambda)(P - t) + \lambda (Q - t)) \varphi_\epsilon (t) \, dt \\
					 & \leq & \int \bigl [(1 - \lambda)F(P - t) + \lambda F(Q - t)\bigr ] \varphi_\epsilon (t) \, dt  \\
					 & = & (1 - \lambda)F_\epsilon(P) + \lambda F_\epsilon(Q) \, .
\end{eqnarray*}
So $F_\epsilon$ is convex.

Now set 
$$
f_j(x) = F_{\epsilon_j} + \delta_j |x|^2 \, .
$$
Certainly $f_j$ is strongly convex because $F_\epsilon$ is convex and $|x|^2$ strongly
convex.  If $\epsilon_j > 0$, $\delta_j > 0$ are chosen appropriately, then we will
have
$$
f_1 \geq f_2 \geq \dots
$$
and $f_j \ra F$ pointwise.  That is the desired conclusion.
\endpf 
\smallskip \\

\begin{proposition} \sl
Let $F: \RR^N \ra \RR$ be a continuous function.  Then $F$ is
convex if and only if, for any $\varphi \in C_c^\infty(\RR^N)$ with
$\varphi \geq 0$, $\int \varphi \, dx = 1$, and any $w = (w_1, w_2, \dots, w_N) \in \RR^N$ it holds that
$$
\int_{\RR^N} F(x) \left [ \sum_{j, k} \frac{\partial^2 \varphi}{\partial x_j \partial x_k} (x) w_j w_k  \right ]\, dx \geq 0 \, .
$$
\end{proposition}
{\bf Proof:}  Assume that $F$ is convex.  In the special case that $F \in C^\infty$, we certainly know that
$$
\sum_{j, k} \frac{\partial^2 F}{\partial x_j \partial x_k} (x) w_j w_k \geq 0 \, .
$$
Hence it follows that
$$
\int_{\RR^N} \sum_{j, k} \frac{\partial^2 F}{\partial x_j \partial x_k} (x) w_j w_k \cdot \varphi(x) \, dx \geq 0 \, .
$$
Now the result follows from integrating by parts twice (the boundary terms vanish since $\varphi$ is
compactly supported).

Now the general case follows by approximating $F$ as in the preceding lemma.

For the converse direction, we again first treat the case when $F \in C^\infty$.
Assume that 
$$
\int_{\RR^N} \sum_{j, k} \frac{\partial^2 \varphi}{\partial x_j \partial x_k} (x) w_j w_k \cdot F(x) \, dx \geq 0 
$$
for all suitable $\varphi$.  Then integration by parts twice gives us the inequality we want.

For general $F$, let $\psi$ be a nonnegative $C_c^\infty$ function, supported in
the unit ball, with integral 1. Set $\psi_\epsilon(x) = \epsilon^{-N} \psi(x/\epsilon)$.  
Define $F_\epsilon(x) = F * \psi_\epsilon(x)$.  Then $F_\epsilon \ra F$ pointwise and
$$
\int_{\RR^N} \sum_{j, k} \frac{\partial^2 \varphi}{\partial x_j \partial x_k} (x) w_j w_k \cdot F(x) \, dx \geq 0 
$$
certainly implies that
$$
\int_{\RR^N} \sum_{j, k} \frac{\partial^2 \varphi}{\partial x_j \partial x_k} (x) w_j w_k \cdot F_\epsilon(x) \, dx \geq 0 \, .
$$
We may integrate by parts twice in this last expression to obtain
$$
\int_{\RR^N} \varphi(x) \sum_{j, k} \frac{\partial^2 F_\epsilon}{\partial x_j \partial x_k} (x) w_j w_k \, dx \geq 0 \, .
$$
It follows that each $F_\epsilon$ is convex.  Thus
$$
F_\epsilon((1 - \lambda)P + \lambda Q) \leq (1 - \lambda)F_\epsilon(P) + \lambda F_\epsilon(Q)
$$
for every $P, Q, \lambda$.  Letting $\epsilon \ra 0^+$ yields that
$$
F((1 - \lambda)P + \lambda Q) \leq (1 - \lambda)F(P) + \lambda F(Q)
$$
hence $F$ is convex.  That completes the proof.
\endpf
\smallskip \\				       

For applications in the next theorem, it is useful to note the following:

\begin{proposition} \sl
Any convex function $f$ is
subharmonic. 
\end{proposition}
{\bf Proof:}  To see this, let $P$ and $P'$ be distinct points in the domain of $f$ and
let $X$ be their midpoint.  Then certainly
$$
2 f(X) \leq f(P) + f(P') \, .
$$
Let $\eta$ be any special orthogonal rotation centered at $X$.  We may write
$$
2 f(X) \leq f(\eta(P)) + f(\eta(P')) \, .
$$
Now integrate out over the special orthogonal group to derive the usual sub-mean-value property
for subharmonic functions.
\endpf
\smallskip \\

The last topic is also treated quite elegantly in Chapter 3 of [HOR].
One may note that the condition that the Hessian be positive semi-definite is stronger
than the condition that the Laplactian be nonnegative.  That gives another proof of
our result.

\begin{theorem} \sl
A domain $\Omega \ss \RR^N$ is convex if and only if it
has a $C^\infty$, strictly convex exhaustion function.
\end{theorem}
{\bf Proof:}  Only the forward direction need be proved (as the
converse direction is contained in the last theorem).

We build the function up iteratively.   We know by the preceding
theorem that there is a continuous exhaustion function $\lambda$.
Let 
$$
\Omega_c = \{x \in \Omega:  \lambda(x) + |x|^2 < c\}
$$
for $c \in \RR$.  Then each $\Omega_c \sss \Omega$ and $c' > c$ implies
that $\Omega_c \sss \Omega'_c$.  Now let $0 \leq \varphi \in C_c^\infty(\RR^N)$ with
$\int \varphi \, dx = 1$, $\varphi$ radial.  We may take $\varphi$ to be supported
in $B(0,1)$.  Let $0 < \epsilon_j < \hbox{dist}(\Omega_{j+1}, \partial \Omega)$.
If $x \in \Omega_{j+1}$, set
$$
\lambda_j(x) = \int_\Omega [\lambda(t) + |t|^2] \epsilon_j^{-N} \varphi ((x - t)/\epsilon_j) \, dV(t) + |x|^2 + 1 \, .
$$
Then each $\lambda_j$ is $C^\infty$ and strictly convex on $\Omega_{j+1}$.   Moreover,
by the previously noted subharmonicity of $\lambda$, we may be sure that
$\lambda_j(x) > \lambda(x) + |x|^2$ on $\overline{\Omega}_j$.

Now let $\chi \in C^\infty(\RR)$ be a convex function with $\chi(t) = 0$ for $t \leq 0$ and 
$\chi'(t), \chi''(t) > 0$ when $t > 0$.  Note that, $\Psi_j(x) \equiv \chi(\lambda_j(x) - (j-1))$ 
is positive and convex on $\Omega_j \setminus \overline{\Omega}_{j-1}$ and is, of course,
$C^\infty$.   Notice now that $\lambda_0 > \lambda$ on $\Omega_0$.  If $a_1$ is large
and positive, then $\lambda'_1 \equiv \lambda_0 + a_1 \Psi_1 > \lambda$ on $\Omega_1$.  
Inductively, if $a_1, a_2, \dots a_{\ell - 1}$ have been chosen, select $a_\ell > 0$ such that
$\lambda'_\ell \equiv \lambda_0 + \sum_{j=1}^\ell a_j \Psi_j >\lambda$ on $\Omega_\ell$.

Since $\Psi_{\ell + k} = 0$ on $\Omega_\ell$, $k > 0$, we see that
$\lambda'_{\ell + k} = \lambda'_{\ell + k'}$ on $\Omega_\ell$ for
any $k, k' > 0$.  So the sequence $\lambda'_\ell$ stabilizes on compacta
and $\lambda' \equiv \lim_{\ell \ra \infty} \lambda'_\ell$ is a $C^\infty$ strictly
convex function that majorizes $\lambda$.  Hence $\lambda'$ is the smooth, strictly convex 
exhaustion function that we seek.
\endpf
\smallskip \\

\begin{corollary} \sl
Let $\Omega \ss \RR^N$ be any convex domain.  Then we may write
$$
\Omega = \bigcup_{j=1}^\infty \Omega_j \, ,
$$
where this is an increasing union and each $\Omega_j$ is strongly convex
with $C^\infty$ boundary.
\end{corollary}
{\bf Proof:}  Let $\lambda$ be a smooth, strictly convex exhaustion function for
$\Omega$.  By Sard's theorem (see [KRP1]), there is a strictly increasing sequence of
values $c_j \ra + \infty$ so that
$$
\Omega_{c_j} = \{x \in \Omega: \lambda(x) < c_j\}
$$
has smooth boundary.  Then of course each $\Omega_{c_j}$ is strongly
convex.  And the $\Omega_{c_j}$ form an increasing sequence of
domains whose union is $\Omega$.
\endpf
\smallskip \\

\section{Other Characterizations of Convexity}

\indent Let $\Omega \ss \RR^N$ be a domain and let ${\cal F}$ be a family
of real-valued functions on $\Omega$ (we do not assume in advance that ${\cal F}$
is closed under any algebraic operations, although often in practice
it will be).  Let $K$ be a compact subset of $\Omega.$  Then the 
{\it convex hull of $K$ in $\Omega$ with respect to ${\cal F}$} is defined to be 
$$
 \widehat{K}_{{\cal F}} \equiv \left \{ x \in \Omega: f(x) \leq 
                                 \sup_{t \in K} f(t) \ \mbox{for all}\  f \in {\cal F} \right \} . 
$$
We sometimes denote this hull by $\widehat{K}$ when the
family ${\cal F}$ is understood or when no confusion is possible.
We say that $\Omega$ is {\it convex} with respect to ${\cal F}$ provided
$\widehat{K}_{\cal F}$ is compact in $\Omega$ whenever $K$ is.  When the
functions in ${\cal F}$ are complex-valued then $|f|$ replaces $f$
in the definition of $\widehat{K}_{\cal F}.$

\begin{proposition}       
Let $\Omega \sss \RR^N$ and let ${\cal F}$ be the family of
real linear functions.  Then $\Omega$ is convex with respect to
${\cal F}$ if and only if $\Omega$ is geometrically convex.
\end{proposition}
{\bf Proof:}  Exercise.  Use the classical definition of convexity
at the beginning of the paper. 
\endpf  \smallskip \\

\begin{proposition}    
Let $\Omega \sss \RR^N$ be any domain.  Let ${\cal F}$ be the family of continuous
functions on $\Omega.$  Then $\Omega$ is convex with respect to ${\cal F}.$
\end{proposition}
{\bf Proof:}  If $K \sss \Omega$ and $x \not \in K$ then the function
$F(t) = 1/(1 + |x - t|)$ is continuous on $\Omega.$  Notice that $f(x) = 1$ and
$|f(k)| < 1$ for all $k \in K.$  Thus $x \not \in \widehat{K}_{\cal F}.$
Therefore $\widehat{K}_{\cal F} = K$ and $\Omega$
is convex with respect to $\cal F.$
\endpf  
\smallskip \\ 

We close this discussion of convexity with a geometric
characterization of the property.  We shall, later in the book,
refer to this as the ``segment characterization''.
First, if $\Omega \ss \RR^N$ is a domain
and $I$ is a closed one-dimensional segment lying in $\Omega$ then
the boundary $\partial I$ is the set consisting of the two endpoints
of $I.$  Now the domain $\Omega$ is convex if and only if whenever
$\{I_j\}_{j=1}^\infty$ is a collection of closed segments in $\Omega$
and $\{\partial I_j\}$ is relatively compact in $\Omega$ then so
is $\{I_j\}.$  This is little more than a restatement of the
classical definition of geometric convexity.  We invite the reader
to supply the details.

In fact the formulation in the last paragraph admits of many variants.  One
of these is the following:  If $\{I_j\}$ is a collection of closed
segmenets in $\Omega$ then
$$
\hbox{dist}\, (\partial I_j, \partial \Omega)
$$
is bounded from 0 if and only if
$$
\hbox{dist}\, (I_j, \partial \Omega)
$$
is bounded from 0.  The following example puts these ideas in perspective.

\begin{example} \rm
Let $\Omega \ss \RR^2$ be 
$$
\Omega = B((0,0),2) \setminus \overline{B((1,0),1)} \, .
$$
Let 
$$
I_j = \{ (- 1/j, t): - 1/2 \leq t \leq 1/2\} \, .
$$
Then it is clear that 
$$
\{\partial I_j\}
$$
is relatively compact in $\Omega$ while 
$$
\{I_j\}
$$
is not.	  And of course $\Omega$ is not convex.
\end{example}

\section{Convexity of Finite Order}

There is a fundamental difference between the domains
$$
B = \{x = (x_1, x_2) \in \RR^2: x_1^2 + x_2^2 < 1\}
$$
and
$$
E = \{x = (x_1, x_2) \in \RR^2: x_1^2 + x_2^4 < 1\} \, .
$$
Both of these domains are convex.
The first of these is strongly convex and the second is not.  More generally,
each of the domains
$$
E_m = \{x = (x_1, x_2) \in \RR^2: x_1^2 + x_2^{2m} < 1\}
$$
is, for $m = 2, 3, \dots$, weakly (not strongly) convex.  Somehow the intuition
is that, as $m$ increases, the domain $E_m$ becomes {\it more} weakly convex.
Put differently, the boundary points $(\pm 1, 0)$ are becoming flatter and
flatter as $m$ increases.

We would like to have a way of quantifying, indeed of measuring, the indicated flatness.
These considerations lead to a new definition.  We first need a bit of terminology.

Let $f$ be a function on an open set $U \ss \RR^N$ and let $P \in \Omega$.
We say that $f$ {\it vanishes to order $k$ at $P$} if any derivative
of $f$, up to and including order $k$, vanishes at $P$.  Thus if $f(P) = 0$
but $\nabla f(P) \ne 0$ then we say that $f$ vanishes to order 0.  If $f(P) = 0$,
$\nabla f(P) = 0$, $\nabla^2 f(P) = 0$, and $\nabla^3 f(P) \ne 0$, then we say 
that $f$ vanishes to order 2.

Let $\Omega$ be a domain and $P \in \partial \Omega$.  Suppose that $\partial \Omega$ is
smooth near $P$.  We say that the tangent plane $T_P(\partial \Omega)$ has order of
contact $k$ with $\partial \Omega$ at $P$ if the defining function $\rho$ for $\Omega$
satisfies
$$
|\rho(x)| \leq C |x - P|^k \qquad \hbox{for all} \ x \in T_P(\partial \Omega) \, ,
$$
and this same inequality does {\it not} hold with $k$ replaced by $k + 1$.

\begin{definition}
Let $\Omega \ss \RR^N$ be a domain and $P \in \partial \Omega$ a point at which
the boundary is at least $C^k$ for $k$ a positive integer.  We say that $P$ is
{\it convex of order $k$} if
\begin{itemize}
\item The point $P$ is convex;
\item The tangent plane to $\partial \Omega$ at $P$ has order of contact $k$ with the boundary at $P$.
\end{itemize}
\end{definition}

\begin{example} \rm
Notice that a point of strong convexity will be convex of order 2.  
The boundary point $(1,0)$ of the domain
$$
E_2k = \{(x_1, x_2) \in \RR^2: x_1^2 + x_2^{2k} < 1\}
$$
is convex of order $2k$.
\end{example}

\begin{proposition} \sl
Let $\Omega \ss \RR^N$ be a bounded domain, and let $P \in \partial \Omega$
be convex of finite order.  Then that order is an even number.
\end{proposition}
{\bf Proof:}  Let $m$ be the order of the point $P$.

We may assume that $P$ is the origin and that the outward
normal direction at $P$ is the $x_1$ direction.  If $\rho$ is a defining
function for $\Omega$ near $P$ then we may use the Taylor expansion about $P$ to write
$$
\rho(x) = 2 x_1 + \varphi(x) \, ,
$$
and $\varphi$ will vanish to order $m$.  If $m$ is odd, then the
domain will not lie on one side of the tangent hyperplane
$$
T_P(\partial \Omega) = \{x: x_1 = 0\} \, .
$$
So $\Omega$ cannot be convex.
\endpf
\smallskip \\

A very important feature of convexity of finite order is its stability.  We formulate
that property as follows:

\begin{proposition} \sl
Let $\Omega \ss \RR^N$ be a smoothly bounded domain and let $P \in \partial \Omega$
be a point that is convex of finite order $m$.  Then points in $\partial \Omega$ that
are sufficiently near $P$ are also convex of finite order {\it at most} $m$.
\end{proposition}
{\bf Proof:}  Let $\Omega = \{x \in \RR^N: \rho(x) < 0\}$, where $\rho$ is a defining
function for $\Omega$.  Then the ``finite order'' condition is given by the
nonvanishing of a derivative of $\rho$ at $P$.  Of course that same derivative
will be nonvanishing at nearby points, and that proves the result.
\endpf  
\smallskip \\

\begin{proposition} \sl
Let $\Omega \ss \RR^N$ be a smoothly bounded domain.  Then there will
be a point $P \in \partial \Omega$ and a neighborhood $U$ of $P$ so
that each point of $U \cap \partial \Omega$ will be convex of order 2
(i.e., strongly convex).
\end{proposition}
{\bf Proof:}  Let $D$ be the diameter of $\Omega$.  We may assume that $\overline{\Omega}$ is distance at
least $10D + 10$ from the origin 0.   Let $P$ be the point of $\partial \Omega$ which is furthest (in the Euclidean
metric) from 0.   Then $P$ is the point that we seek.  Refer to Figure 1.

Let $L$ be the distance of 0 to $P$.  Then we see that the sphere with center 0 and radius $L$
externally osculates $\partial \Omega$ at $P$.  Of course the sphere is strongly convex
at the point of contact.  Hence so is $\partial \Omega$.  By the continuity of second
derivatives of the defining function for $\Omega$, the same property holds for
nearby points in the boundary.  That completes the proof.
\endpf
\smallskip \\

\begin{example} \rm
Consider the domain
$$
\Omega = \{(x_1, x_2, x_3) \in \RR^3: x_1^2 + x_2^4 + x_3^4 < 1\} \, .
$$
The boundary points of the form $(0,a,b)$ are convex of order 4.  All
others are convex of order 2 (i.e., strongly convex).
\end{example}

It is straightforward to check that Euclidean isometries preserve convexity, preserve
strong convexity, and preserve convexity of finite order.  Diffeomorphisms do not.
In fact we have:

\begin{proposition} \sl
Let $\Omega_1$, $\Omega_2$ be smoothly bounded domains in $\RR^N$,
let $P_1 \in \partial \Omega_1$ and $P_2 \in \partial \Omega_2$.   Let
$\Phi$ be a diffeomorphism from $\overline{\Omega_1}$ to $\overline{\Omega_2}$
and assume that $\Phi(P_1) = P_2$.  Further suppose that
the Jacobian matrix of $\Phi$ at $P_1$ is an orthogonal linear
mapping.  Then we have:
\begin{itemize}
\item If $P_1$ is a convex boundary point then $P_2$ is a convex boundary point;
\item If $P_1$ is a strongly convex boundary point then $P_2$ is a strongly convex boundary point;
\item If $P_1$ is a boundary point that is convex of order $2k$ then $P_2$ is a boundary point that is 
convex of order $2k$.
\end{itemize}
\end{proposition}
{\bf Proof:}  We consider the first assertion.  Let $\rho$ be a 
defining function for $\Omega_1$. Then $\rho \circ \Phi^{-1}$ 
will be a defining function for $\Omega_2$.  Of course we know
that the Hessian of $\rho$ at $P_1$ is positive semi-definite.
It is straightforward to calculate the Hessian of $\rho' \equiv \rho \circ \Phi^{-1}$
and see that it is just the Hessian of $\rho$ composed with $\Phi$ applied to
the vectors transformed under $\Phi^{-1}$.  So of course $\rho'$ will
have positive semi-definite Hessian.

The other two results are verified using the same calculation.
\endpf
\smallskip \\

\begin{proposition} \sl
Let $\Omega$ be a smoothly bounded domain in $\RR^N$.  Let ${\cal L}$ be an 
invertible linear map on $\RR^N$.  Define $\Omega' = {\cal L}(\Omega)$.
Then 
\begin{itemize}
\item Each convex boundary point of $\Omega$ is mapped to a convex boundary
point of $\Omega'$.
\item Each strongly convex boundary point of $\Omega$ is mapped to a strongly convex boundary
point of $\Omega'$.
\item Each boundary point of $\Omega$ that is convex of order $2k$ is mapped to a boundary
point of $\Omega'$ that is convex of order $2k$.
\end{itemize}
\end{proposition}
{\bf Proof:}  Obvious.
\endpf 
\smallskip \\

Maps which are not invertible tend to decrease the order of a convex point.
An example will illustrate this idea:

\begin{example} \rm
Let $\Omega = \{(x_1, x_2) \in \RR^2: x_1^2 + x_2^2 < 1\}$ be the unit ball
and $\Omega' = \{(x_1, x_2) \in \RR^2: x_1^4 + x_2^4 < 1\}$.  We see
that
$$
\Phi(x_1, x_2) = (x_1^2, x_2^2)
$$
maps $\Omega'$ onto $\Omega$.  And we see that $\Omega$ is strongly convex (i.e., convex
of order 2 at each boundary point) while $\Omega'$ has boundary points that are convex
of order 4.  The points of order 4 are mapped by $\Phi$ to points of order 2.
\end{example}

\section{Extreme Points}

A point $P \in \partial \Omega$ is called an {\it extreme point} if,
whenever $a, b \in \partial\Omega$ and $P = (1- \lambda)a + \lambda b$
for some $0 \leq \lambda \leq 1$ then $a = b = P$.

It is easy to see that, on a convex domain, a point of strong
convexity must be extreme, and a point that is convex of order
$2k$ must be extreme. But convex points in general are {\it
not} extreme.

\begin{example} \rm
Let 
$$								(
\Omega = \{(x_1, x_2) \in \RR^2: |x_1| < 1, |x_2| < 1\} \, .
$$
Then $\Omega$ is clearly convex.  But any boundary
point with $x_1$, $x_2$ not both 1 is not extreme.

For example, consider the boundary point $(1/2, 1)$.  Then
$$
(1/2, 1) = \frac{1}{2} (1/4, 1) + \frac{1}{2} (3/4, 1) \, .
$$
\end{example}

\begin{example}  \rm
Let $\Omega \ss \RR^2$ be the domain with boundary
consisting of
\begin{itemize}
\item The segments from $(-3/4,1)$ to $(3/4,1)$,
from $(1,3/4)$ to $(1, -3/4)$, from $(3/4,-1)$ to $(-3/4,-1)$,
and from $(-1, 3/4)$ to $(-1, -3/4)$.
\item The four circular arcs
$$
(x + 3/4)^2 + (y - 3/4)^2 = \frac{1}{16} \ , \quad y \geq 0 \, , \ x \leq 0 \, ;
$$
$$
(x - 3/4)^2 + (y - 3/4)^2 = \frac{1}{16} \ , \quad y \geq 0 \, , \ x \geq 0 \, ;
$$
$$
(x - 3/4)^2 + (y + 3/4)^2 = \frac{1}{16} \ , \quad y \leq 0 \, , \ x \geq 0 \, ;
$$
$$
(x + 3/4)^2 + (y + 3/4)^2 = \frac{1}{16} \ , \quad y \leq 0 \, , \ x \leq 0 \, ;
$$
\end{itemize}
Then any point on any of the circular arcs is extreme.  But no other boundary
point is extreme.  Note, however, that the extreme points $(-3/4,1)$, $(3/4,1)$,
$(1, -3/4)$, $(1, 3/4)$, $(-3/4,-1)$, $(3/4, -1)$, $(-1, -3,4)$, and $(-1,3/4)$ are {\it not}
convex of finite order.  
\end{example}

\section{Support Functions}

Let $\Omega \ss \RR^N$ be a bounded, convex domain with $C^2$ boundary.
If $P \in \partial \Omega$ then let $T_P(\partial \Omega)$ be the tangent
hyperplane to $\partial \Omega$ at $P$.  We may take the outward unit normal
at $P$ to be the positive $x_1$ direction.   Then the function
$$
L(x) = x_1
$$
is a linear function that is negative on $\Omega$ and positive on the
other side of $T_P(\Omega)$.  The function $L$ is called a {\it support function}
for $\Omega$ at $P$.  Note that if we take the supremum of all support functions
for all $P \in \partial \Omega$ then we obtain a defining function for $\Omega$.

The support function of course takes the value 0 at $P$.  It may take the
value 0 at other boundary points---for instance in the case of the domain
$\{(x_1, x_2): |x_1| < 1, |x_2|< 1\}$.   But if $\Omega$ is convex and $P \in \partial \Omega$
is a point of convexity of finite order $2k$ then the support function will vanish on
$\partial \Omega$ only at the point $P$.   The same assertion holds when $P$ 
is an extreme point of the boundary.

\section{Bumping}

One of the features that distinguishes a convex point of finite order from a convex point
of infinite order is stability.  The next example illustrates the point.

\begin{example} \rm
Let 
$$
\Omega = \{(x,y) \in \RR^2: |x| < 1, |y| < 1\} \, .
$$
Let $P$ be the boundary point $(1/2, 1)$.  Let $U$ be a small
open disc about $P$.  Then there is no open domain $\widehat{\Omega}$
such that
\begin{enumerate}
\item[{\bf (a)}]  $\widehat{\Omega} \supseteq \Omega$ and $\widehat{\Omega} \ni P$;
\item[{\bf (b)}]  $\widehat{\Omega} \setminus \Omega \ss U$;
\item[{\bf (b)}]  $\widehat{\Omega}$ is convex.
\end{enumerate}
To see this assertion, assume not.  Let $x$ be a point of $\widehat{\Omega} \setminus \Omega$.
Let $y$ be the point $(0.9, 0.9) \in \Omega$.  Then the segment connecting $x$ with $y$ will {\it not}
lie in $\widehat{\Omega}$.  
\end{example}

The example shows that a flat point in the boundary of a convex domain cannot be perturbed
while preserving convexity.  But a point of finite order {\it can} be perturbed:

\begin{proposition} \sl Let $\Omega \ss \RR^N$ be a 
bounded, convex domain with $C^k$ boundary. Let $P \in \partial \Omega$ be a
convex point of finite order $m$. Write $\Omega = \{x \in
\RR^N: \rho(x) < 0\}$. Let $\epsilon > 0$. Then there is a
perturbed domain $\widehat{\Omega} = \{x \in
\RR^N: \widehat{\rho}(x) < 0\}$ with $C^k$ boundary such that 
\begin{enumerate}
\item[{\bf (a)}] $\widehat{\Omega} \supseteq \Omega$;
\item[{\bf (b)}] $\widehat{\Omega} \ni P$; 
\item[{\bf (c)}]
$\partial\widehat{\Omega} \setminus \overline{\Omega}$
consists of points of finite order $m$; 
\item[{\bf (d)}] The Hausdorff distance of $\partial \widehat{\Omega}$ and
$\partial\Omega$ is less than $\epsilon$. 
\end{enumerate}
\end{proposition} 

\noindent Before we begin the proof, we provide a useful technical
lemma:

\begin{lemma} \sl  Let $a$ be a fixed, positive number.   Let $\alpha_0, \alpha_1, \dots,
\alpha_k$ and $\beta_0, \beta_1, \dots, \beta_k$ and $\gamma_0$ be real numbers.  Then
there is a concave-down polynomial polynomial function $y = p(x)$ so that
\begin{itemize}
\item $p(0) = \gamma_0$;
\item $p(-a) = \alpha_0$, $p(a) = \beta_0$;
\item $p^{(j)}(-a) = \alpha_j$ for $j = 1, \dots, k$;
\item $p^{(j)}(a) = \beta_j$ for $j = 1, \dots, k$.
\end{itemize}
Here the exponents in parentheses are derivatives.
\end{lemma}
{\bf Proof of the Lemma:}  Define
$$
g_j(x) = x^{2j}
$$
and let $h_j^{\theta_j}$ be the function obtained from $g_j$ by rotating
the coordinates $(x,y)$ through an angle of $\theta_j$.  Define
$$
p(x) = c_0 - (c_1)^2 h_1^{\theta_1}(x) - (c_2)^2 h_2^{\theta_2}(x) 
	   - \cdots - (c_p)^2 h_p^{\theta_p}(x) \, ,
$$
some positive integer $p$.  If $p$ is large enough (at least $k+1$),
then there will be more free parameters in the definition of
$p$ than there are constants $\alpha_j$, $\beta_j$, and $\gamma_0$.  So we
may solve for the $c_j$ and $\theta_j$ and thereby define $p$.
\endpf 
\smallskip \\

\noindent {\bf Proof of the Proposition:} First let us consider the case
$N = 2$. Fix $P \in \partial \Omega$ as given in the statement
of the proposition.  We may assume without loss of generality that $P$ is the
origin and the tangent line to $\partial \Omega$ at $P$ is the
$x$-axis.  We may further assume that $\Omega$ is so oriented that
the boundary $\partial \Omega$ of $\Omega$ near $P$ is the graph of
a concave-down function $\varphi$.

Let $\delta > 0$ be small and let $x$ and
$y$ be the two boundary points that are horizontally distance $\delta$ from
$P$ (situated, respectively, to the left and to the right of $P$). If $\delta$ is sufficiently small, then the angle
between the tangent lines at $x$ and at $y$ will be less than
$\pi/6$.

Now we think of $P = (0,0)$, $\gamma_0 = \epsilon > 0$, of $x = (-a, \alpha_0)$, and
of $y = (a, \beta_0)$.  Further, we set
$$
\alpha_j = \varphi^{(j)}(-a) \ , \quad j = 1, \dots, k
$$
and
$$
\beta_j = \varphi^{(j)}(a) \ , \quad j = 1, \dots, k \, .
$$
Then we may apply the lemma to obtain a concave-down polynomial $p$ 
which agrees with $\varphi$ to order $k$ at the points of
contact $x$ and $y$.  

Thus the domain $\widehat{\Omega}$ which has boundary given
by $y = p(x)$ for $x \in [-a,a]$ and boundary coinciding with 
$\partial \Omega$ elsewhere (we are simply replacing the
portion of $\partial \Omega$ which lies between $x$ and $y$ by
the graph of $p$) will be a convex domain that bumps $\Omega$ 
{\it provided that} the degree of $p$ does not exceed the finite
order of convexity $m$ of $\partial \Omega$ near $P$.   When the
degree of $p$ exceeds $m$, then the graph $y = p(x)$ may intersect
$\partial \Omega$ between $x$ and $y$, and therefore not provide
a geometrically valid bump.

For higher dimensions, we proceed by slicing.  Let $P \in \partial \Omega$ 
be of finite order $m$.  Let $T_P(\partial \Omega)$ be the tangent hyperplane to $\partial \Omega$
at $P$ as usual.  If {\bf v} is a unit vector in $T_P(\partial \Omega)$ and $\nu_P$ the unit
outward normal vector to $\partial \Omega$ at $P$, then consider the 2-dimensional plane
${\cal P}_{\bf v}$ spanned by {\bf v} and $\nu_P$.  Then $\Omega_{\bf v} \equiv {\cal P}_{\bf v} \cap \Omega$
is a 2-dimensional convex domain which is convex of order $m$ at $P \in \partial \Omega_{\bf v}$.
We may apply the two-dimensional perturbation result to this domain.  We do so for each unit tangent vector
${\bf v} \in T_P(\partial \Omega)$, noting that the construction varies smoothly with the data
vector {\bf v}.  The result is a smooth, perturbed domain $\widehat{\Omega}$ as desired.
\endpf 
\smallskip \\
	     
It is worth noting that the proof shows that, when we bump a piece of boundary that 
is convex of order $m$, then we may take the bump to be convex of order 2 or 4 or any
degree up to and including $m$ (which of course is even)w.

It is fortunate that the matter of bumping may be treated more or less heuristically in the present
context.  In several complex variables, 
bumping is a more profound and considerably
more complicated matter (see, for instance [BHS]).  

\section{Concluding Remarks}

We have attempted here to provide the analytic tools so that convexity
can be used in works of geometric analysis.  There are many other
byways to be explored in this vein, and we hope to treat them at another time.

\vspace*{.48in}

\noindent {\Large \sc References}
\smallskip  \\

\begin{enumerate}

\item[{\bf [BHS]}] Gautam Bharali and Berit Stens\o nes,
Plurisubharmonic polynomials and bumping, {\it Math.\ Z.}
261(2009), 39--63.

\item[{\bf [BOF]}] T. Bonneson and W. Fenchel, {\it Theorie der
konvexen K\"{o}rper}, Springer-Verlag, Berlin, 1934.

\item[{\bf [FEN]}] W. Fenchel, Convexity through the ages, {\it
Convexity and its Applications}, Birkh\"{a}user, Basel, 1983,
120--130.

\item[{\bf [HOR]}]  L. H\"{o}rmander,  {\it Notions of Convexity}, Birkh\"{a}user Publishing,
Boston, MA, 1994.

\item[{\bf [KRA1]}]  S. G. Krantz, {\it Function Theory of Several Complex
Variables}, $2^{\rm nd}$ ed., American Mathematical Society, Providence,
RI, 2001.

\item[{\bf [KRP1]}]  S. G. Krantz and H. R. Parks, {\it The Geometry of Domains
in Space}, Birkh\"{a}user Publishing, Boston, MA, 1999.

\item[{\bf [KRP2]}]  S. G. Krantz and H. R. Parks, {\it The Implicit Function
Theorem}, Birkh\"{a}user, Boston, 2002.
	  
\item[{\bf [LAY]}] S. R. Lay, {\it Convex Sets and Their
Applications}, John Wiley and Sons, New York, 1982.

\item[{\bf [LEM]}]  L. Lempert, La metrique \mbox{K}obayashi et las representation des domains
sur la boule, {\em Bull. Soc. Math. France} 109(1981), 427-474.

\item[{\bf [ONE]}]  B. O'Neill, {\it Elementary Differential
Geometry,} Academic Press, New York, 1966.

\item[{\bf [VAL]}] F. A. Valentine, {\it Convex Sets},
McGraw-Hill, New York, 1964.

\item[{\bf [VLA]}] V. Vladimirov, {\it Methods of the Theory of
Functions of Several Complex Variables}, MIT Press, Cambridge,
1966.				     

\item[{\bf [ZYG]}]  A. Zygmund, {\it Trigonometric Series},
Cambridge University Press, Cambridge, UK, 1968.

\end{enumerate}
\vspace*{.4in}

\begin{quote}
Department of Mathematics \\
Washington University in St.\ Louis \\
St.\ Louis, Missouri 63130 \ \ USA  \\
{\tt sk@math.wustl.edu}  \\
\end{quote}

\end{document}